\documentclass[10pt]{article}
\usepackage{amssymb} 
\usepackage{amsmath} 
\usepackage{latexsym}
\newtheorem{Proof}{Proof.}
 
 \newenvironment{proof}{\begin{Proof}\rm}{\hfill $\Box$ \end{Proof}}
\newtheorem{theorem}{Theorem}

\newtheorem{ax}{A}

\begin{document}
\title{A Problem in Pythagorean Arithmetic}
\author{Victor Pambuccian}
\date{}
\maketitle

\begin{abstract}
Problem 2 at the 56th International Mathematical Olympiad (2015)
asks for all triples $(a,b,c)$ of positive integers for which
$ab-c$, $bc-a$, and $ca-b$ are all powers of $2$. We show that this
problem requires only a primitive form of arithmetic, going back to
the Pythagoreans, which is the arithmetic of the even and the odd.
\end{abstract}

\section{Introduction}

Problem 2 at the 56th International Mathematical Olympiad (2015),
proposed by Du\v{s}an Djuki\'{c}, asked contestants to find all
triples $(a,b,c)$ of positive integers for which $ab-c$, $bc-a$, and
$ca-b$ are all powers of $2$. Here a ``power of $2$" is understood
to be $2^n$ with $n$ a non-negative integer.

As is well known, problems at the IMO should be solvable with {\em
elementary} means, and our aim is to find out just how elementary a
formal theory is needed to solve Problem 2. Since it speaks about
positive integers and the operations of addition and multiplication,
an axiom system for a theory in which it holds will need to contain
the binary operations $+$ and $\cdot$, the binary relation $<$, as
well as the constants $0$ (so that we can express that all numbers
we deal with are non-negative) and $1$ (so that we can express the
fact that the successor of a number $n$ in the order determined by
$<$ is $n+1$).

\section{The axiom system for PA$^-$ and its extensions}

Thus we need axioms for the usual rules for addition $+$ and
multiplication $\cdot$, for $1$ and  $0$, that is:

\begin{ax} \label{ass+}
$(x+y)+z = x+(y+z)$
\end{ax}
\begin{ax} \label{com+}
$x+y =  y+x$
\end{ax}
\begin{ax} \label{ass.}
$(x\cdot y)\cdot z = x\cdot (y\cdot z)$
\end{ax}
\begin{ax} \label{com.}
 $x\cdot y = y\cdot x$
\end{ax}
\begin{ax} \label{dis}
$x\cdot (y+z) = x\cdot y+x\cdot z$
\end{ax}
\begin{ax} \label{zer}
$x+0=x \wedge x\cdot 0 = 0$
\end{ax}
\begin{ax} \label{one}
$x\cdot 1 = x$
\end{ax}

We also need axioms for inequality $<$, and a binary operation $-$,
so that we can express the difference between two numbers if the
result is positive.  These are

\begin{ax} \label{tra}
 $(x<y \wedge y<z) \rightarrow x<z$
 \end{ax}
\begin{ax} \label{ire}
$\neg x<x$
 \end{ax}
\begin{ax} \label{tot}
$x<y \vee x=y \vee y<x$
\end{ax}
\begin{ax} \label{add}
$x<y \rightarrow x+z < y+z$
\end{ax}
\begin{ax} \label{mul}
$(0<z \wedge x<y) \rightarrow x\cdot z < y\cdot z$
\end{ax}
\begin{ax} \label{min}
$x<y \rightarrow x+(y-x) =y$
\end{ax}
\begin{ax} \label{disc}
$0<1 \wedge (x>0 \rightarrow (x>1 \vee x=1))$
\end{ax}
\begin{ax} \label{pos}
$x>0 \vee x=0$
\end{ax}

A\ref{ass+}-A\ref{pos} represents an axiom system for what is
referred to as PA$^-$ in \cite[pp.\ 16-18]{K91}. Models of PA$^-$
consist of {\em numerals}, i.e., $\overline{n}=(1+(1+\ldots +1))$,
with $1$ ocuurring $n$ times, and, possibly, of {\em nonstandard}
elements, which are greater than all numerals.

By referring to powers of $2$, our problem seems to require more,
for we do not have the exponential function in our vocabulary. It
turns out that we do we need it, for we can express the fact that
$a$ is a power of $2$ simply by defining the unary predicate $PT$
which stipulates that a positive number $a$ is a power of $2$ if and
only if all its divisors, except $1$, are even:

\begin{equation} \label{1}
PT(n):\Leftrightarrow n>0 \wedge (\forall d)\, (d|n \wedge d>1
\rightarrow \overline{2}|d).
\end{equation}

This definition certainly corresponds to our intuitions regarding
powers of $2$, but it must not satisfy properties we find to be
intrinsic to the notion of ``power of $2$", which can be formalized
as follows:

\begin{eqnarray}
PT(a)\wedge PT(b)\rightarrow PT(ab) \label{2ab}\\ PT(a)\wedge
PT(b)\wedge a<b \rightarrow a|b\label{a|b}\\  PT(a)\wedge a<b\wedge
b<\overline{2}\cdot a\rightarrow \neg PT(b)\label{ab2a}
\end{eqnarray}

This is perhaps not so surprising if one thinks that PA$^{-}$ is a
very weak theory, in which one cannot even show that among two
consecutive numbers one is even and the other none is odd. In fact,
for any natural number $n$, there may be sequences of $n$
consecutive numbers, none of which is odd and none of which is even.
For the positive cone of ${\mathbb Z}[X]$ (here ${\mathbb Z}[X]$ is
ordered by $\sum_{i=0}^n c_iX^i>0$ if and only if $c_n>0$) is a
model of PA$^{-}$, and the sequence $X+1, \ldots, X+\overline{n}$
has no even element and no odd element. Yet none of
(\ref{2ab})-(\ref{ab2a}) holds in PA$^-$+A\ref{oe} either, where
A\ref{oe} is the axiom expressed in a language enriched with the
unary operation symbol $\left[\frac{\cdot}{2}\right]$, stating that
every number is odd or even:

\begin{ax} \label{oe}
$x=2\left[\frac{x}{2}\right]\vee x=2\left[\frac{x}{2}\right]+1$.
\end{ax}

To see this, denote by  $K_{D}[X]$ the ring of polynomials in $X$
with free term in $D$ and with all other coefficients in $K$,
ordered by $\sum_{i=0}^n c_iX^i>0$ if and only if $c_n>0$ (here
$c_0\in D$, and $c_i\in K$ for all $1\leq i\leq n$, with $c_n\neq
0$), and  denote by ${\mathcal C}(K_{D}[X])$ the positive cone of
$K_{D}[X]$. Let ${\mathbb Z}_{\frac{1}{2}}$ stand for the ring of
dyadic numbers, i.e., all rational numbers of the form
$\frac{m}{2^n}$, with $m, n\in {\mathbb Z}$ and $n\geq 0$, and let
$R={\mathbb Z}_{\frac{1}{2}}[\sqrt{3}]$ stand for the ring whose
elements are of the form $a+b\sqrt{3}$, with $a, b \in {\mathbb
Z}_{\frac{1}{2}}$. Then ${\mathcal C}(R_{{\mathbb Z}}[X])$ with
$\left[\frac{\sum_{i=1}^n a_iX^i +a_0}{2}\right]=\sum_{i=1}^n
\frac{a_i}{2}X^i + \left[\frac{a_0}{2}\right]$ is a model of PA$^-+
A\ref{oe}$. However, given that $PT(\sqrt{3}X)$, but $\neg
PT(3X^2)$, (\ref{2ab}) does not hold, and given that $PT(X)$,
$PT(\sqrt{3}X)$, $X<\sqrt{3}X$, yet $X\nmid \sqrt{3}X$, (\ref{a|b})
does not hold either, and the fact that $X<\sqrt{3}X<2X$, with $X$,
$\sqrt{3}X$, and $2X$ powers of $2$, shows that (\ref{ab2a}) does
not hold.

What PA$^-+ A\ref{oe}$ lacks is an axiom stating that every fraction
can be brought into a form in which  numerator and denominator are
not both even. It is an axiom needed for the proof based on
considerations of parity of the fact that $\sqrt{2}$ is irrational.
This was, apparently, the oldest form of number theory, as practiced
by the Pythagoreans, about which Aristotle tells us in his {\em
Metaphysics}, 986a, that

\begin{quote}
``Evidently, then, these thinkers also consider that number is the
principle both as matter for things and as forming both their
modifications and their permanent states, and hold that the elements
of number are   the  even and the odd" (translated by W. D. Ross)
\end{quote}

For more on the arithmetic of the even and the odd, see \cite{pamb}.
To state the axiom, we need three more binary operations, $\kappa$,
$\mu$, and $\nu$ (so the language in which our {\em Pythagorean
Arithmetic} is expressed consists of $0$, $1$, $+$, $\cdot$, $<$,
$\left[\frac{\cdot}{2}\right]$, $\kappa$, $\mu$, $\nu$):

\begin{ax} \label{rf}
$m= \kappa(m,n)\cdot \mu(m,n)\wedge n= \kappa(m,n)\cdot
\nu(m,n)$\\
\hspace*{10mm} $\wedge (\mu(m,n)=2\left[\frac{\mu(m,n)}{2}\right] +1
\vee \nu(m,n)=2\left[\frac{\nu(m,n)}{2}\right] +1)$
\end{ax}

Notice that A\ref{oe} becomes superfluous in the presence of
A\ref{rf}, as it follows by  applying A\ref{rf} with $m=2$, and
noticing that, in PA$^-$, if $a\cdot b=2$, then $a=1$ or $b=1$. {\em
Pythagorean Arithmetic} can  thus be axiomatized by
$\{$A\ref{ass+}-A\ref{pos}, A\ref{rf}$\}$.

Throughout the paper, we will use the symbols $\leq$ and $\geq$ with
their usual meanings. All of (\ref{2ab})-(\ref{ab2a}) hold in
Pythagorean Arithmetic. To see this, notice first that cancelation
holds, i.e., satisfies the following

\begin{equation} \label{l-}
a+x=a+y \rightarrow x=y
\end{equation}

\begin{proof}
Suppose $a+x=a+y$. By A\ref{tot}, one of $x<y$, $x=y$, or $y<x$ must
hold. Suppose $x=y$ does not hold. Given the symmetry in $x$ and $y$
of our hypothesis, we may assume, w.\ l.\ o.\ g.\ that $x<y$. Then,
by A\ref{add}, we have $a+x<a+y$ as well, thus $a+y<a+y$, which
contradicts A\ref{ire}.
\end{proof}

Cancelation is allowed, i.e.,

\begin{equation} \label{canc}
a\neq 0\wedge a\cdot x=a\cdot y \rightarrow x=y
\end{equation}

\begin{proof}
By (\ref{tot}), we have $x<y$ or $x=y$, or $y<x$. If $x=y$ does not
hold, then one of $x<y$ or $y<x$ must hold. Suppose $x<y$. By
A\ref{mul}, we have $a\cdot x< a\cdot y$, contradicting our
hypothesis. Same contradiction by assuming $y<x$.
\end{proof}

Distributivity of multiplication holds over subtraction as well,
i.e.,
\begin{equation} \label{dis-}
b<a \rightarrow c\cdot a - c\cdot b = c\cdot (a - b)
\end{equation}

\begin{proof}
By A\ref{min}, $b+ (a-b)=a$, thus, by A\ref{dis}, $c\cdot a=
c\cdot(b+ (a-b))= c\cdot b + c\cdot (a-b)$, and, since $c\cdot b +
(c\cdot a - c\cdot b) = c\cdot a$, we must, by (\ref{l-}), $c\cdot
(a-b)= c\cdot a - c\cdot b$.
\end{proof}

Also, odd numbers are never even, i.e.,

\begin{equation} \label{loe}
\overline{2}\cdot n+1\neq \overline{2}\cdot m
\end{equation}

\begin{proof}
Suppose $\overline{2}\cdot n+1= \overline{2}\cdot m$. By A\ref{disc}
and A\ref{add} $\overline{2}\cdot n<\overline{2}\cdot n+1$, thus,
$\overline{2}\cdot n<\overline{2}\cdot m$, so, by A\ref{min},
$\overline{2}\cdot n+ (\overline{2}\cdot m-\overline{2}\cdot n)
=\overline{2}\cdot m$. Thus, $\overline{2}\cdot n+
(\overline{2}\cdot m-\overline{2}\cdot n)=\overline{2}\cdot n+1$,
and thus, by (\ref{l-}), $\overline{2}\cdot m-\overline{2}\cdot
n=1$, i.e., by (\ref{dis-}), $\overline{2}\cdot (m-n)=1$. Since
$m-n>0$, we have, by A\ref{disc}, $m-n >1$ or $m-n=1$. Thus, by
A\ref{one}, $\overline{2}\cdot (m-n) > \overline{2}$ or
$\overline{2}\cdot (m-n) = \overline{2}$, i.e., $1 > \overline{2}$
or $1 = \overline{2}$, none of which can hold, for, by A\ref{disc}
and A\ref{add}, $0<1$ and $1<1+1$.
\end{proof}

We also have:

\begin{equation} \label{div2}
\overline{2}\cdot m+1 | a\cdot b \wedge PT(a) \rightarrow
\overline{2}\cdot m+1 | b.
\end{equation}

\begin{proof}
Since $\overline{2}\cdot m+1 | a\cdot b$, there must be a $c$ such
that $(\overline{2}\cdot m+1)\cdot c=a\cdot b$. By A\ref{rf} with
$c$ instead of $m$ and $a$ instead of $n$ we get that
$c=\kappa(c,a)\cdot \mu(c,a)$ and $a==\kappa(c,a)\cdot \nu(c,a)$,
with at least one of $\mu(c,a)$ and $\nu(c,a)$ odd. Plugging in to
$(\overline{2}\cdot m+1)\cdot c=ab$ and canceling $\kappa(a,c)$, we
get $(2m+1)\cdot \mu(c,a)=\nu(c,a)\cdot b$. Now $\nu(c,a)$ must be
odd, for, if it were even, $(\overline{2}\cdot m+1)\cdot \mu(c,a)$
would have to be even as well, forcing $\mu(c,a)$ to be even (it has
to be even or odd, since A\ref{oe} holds, and, if it were odd,
$(2m+1)\mu(c,a)$ would be odd, a contradiction, for a number cannot
be both odd and even, by (\ref{loe})), but one of $\nu(c,a)$ and
$\mu(c,a)$ must be odd. Since $\nu(a,c)$ is odd, $\nu(c,a) | a$ and
$PT(a)$, we must have $\nu(c,a)=1$, so we have $(\overline{2}\cdot
m+1)\cdot \mu(c,a)=b$, so $\overline{2}\cdot m+1 | b$.
\end{proof}

We can now show that (\ref{2ab})-(\ref{ab2a}) hold in Pythagorean
Arithemtic. Suppose $PT(a)$ and $PT(b)$ and let $d | ab$ with $d>1$.
If $d$ were odd, then, by (\ref{div2}), bearing in mind that
$PT(a)$, we would have $d | b$, but that would contradict the fact
that $PT(b)$. This proves (\ref{2ab}). Suppose now $a<b$, $PT(a)$
and $PT(b)$. By A\ref{rf} we have $a=\kappa(a,b)\cdot \mu(a,b)$ and
$b=\kappa(a,b)\cdot \nu(a,b)$. Since $a$ and $b$ cannot have odd
divisors greater than $1$, and one of $\mu(a,b)$  and $\nu(a,b)$ has
to be odd, the odd one has to be $1$ (both cannot be $1$, for else
$a=b$). Since we cannot have $\nu(a,b)=1$, as that would entail
$b<a$ or $b=a$, we must have $\mu(a,b)=1$, and thus $a | b$, proving
(\ref{a|b}). Suppose now $a<b$, $b<\overline{2}\cdot a$, and
$PT(a)$.  By A\ref{rf}, we have $a=\kappa(a,b)\cdot \mu(a,b)$ and
$b=\kappa(a,b)\cdot \nu(a,b)$. Given that $a$ can have no odd
divisor except for $1$, $\mu(a,b)$ is either even or $1$. If it were
$1$, then $b=a\cdot \nu(a,b)$, and thus $1<\nu(a,b)<\overline{2}$,
contradicting A\ref{disc}, which asks for $\nu(a,b)-1$ to be $1$ or
$>1$, i.e. $\nu(a,b)=\overline{2}$ or $\nu(a,b)>\overline{2}$, a
contradiction. Thus  $\mu(a,b)$ is  even, so  $\nu(a,b)$ must be
odd. It cannot be $1$, for else we would have $b\leq a$, so
 $\nu(a,b)$ is an odd number greater than $1$. Thus $\neg PT(b)$,
 proving (\ref{ab2a}).

\section{Problem 2 holds in Pythagorean Arithmetic}

To turn Problem 2 into a statement that can be proved inside
Pythagorean Arithmetic, we need to express it not as a question but
rather as a solved problem, one that states what that solutions are
and implicitly that there are no other solutions. In this form, its
statement is --- with $S=\{(2,2,2), (3,2,2), \&, (11,6,2), \&,
(7,5,3), \&\}$, where by $(x,y,z), \&$ we have denoted the sequence
of all triples obtained by permuting $x,y,$ and $z$ ---

\begin{eqnarray} \label{p2}
& & a\cdot b>c \wedge b\cdot c>a \wedge  c\cdot a > b\wedge
PT(a\cdot
b-c)\wedge PT(b\cdot c-a)\wedge PT(c\cdot a - b)\nonumber\\
\hspace*{9mm} & & \rightarrow \bigvee_{(i,j,k)\in S} a=i\wedge b=j
\wedge c=k
\end{eqnarray}

\begin{theorem}
The statement (\ref{p2}) can be proved using only the axioms
$\{$A\ref{ass+}-A\ref{pos}, A\ref{rf}$\}$, i.e., inside
Pytha\-gorean Arithmetic.
\end{theorem}

\begin{proof}
First, notice that each of $a, b,$ and $c$ has to be greater than
$1$. That none can be $0$ is plain, for if, say $a=0$, then $a\cdot
b>c$ could not hold, given A\ref{pos}. None of them can be $1$
either, for if, say, $a=1$, then we would have $b>c$ and $c>b$,
which, after applying A\ref{tra}, would contradict A\ref{ire}.
Suppose now that two of $a, b,$ and $c$ were equal, say, $a=b$. Then
we would have $PT(a^2-c)$ and $PT(a\cdot(c-1)$. The latter implies
$PT(a)$ and $PT(c-1)$, and, since $a>1$, $PT(a)$ implies that $a$ is
even. If $c>2$, then $c-1>1$, and thus $PT(c-1)$ would imply that
$c-1$ is even, i.e., $c$ is odd. But then $a^2-c$ would have to be
odd, and since we have $P(a^2-c)$, we would need to have $a^2-c=1$,
i.e., $a^2=c+1$. Since $PT(a)$, we also have, by (\ref{2ab}),
$PT(a^2)$, so $PT(c+1)$ as well. Given that their difference is
$\overline{2}$, both $c-1$ and $c+1$, which have to be even, as
$c>\overline{2}$, cannot be multiples of $\overline{4}$. Since both
have only even divisors, one of them must be $\overline{2}$. Since
$c+1>\overline{3}$, we must have $c-1=\overline{2}$, so
$c=\overline{3}$, and thus, given $a^2=c+1$, $a=\overline{2}$. So
$(2,2,3)$ is the only solution with $a=b$ and $c>2$. If $c=2$, then
$PT(a^2-c)$ and $PT(a)$ imply that $4\nmid a$, so that $a=2$. Thus
$(2,2,2)$ is the only solution with $a=b$ and $c=2$.

Given the symmetry in $a,b,c$ of the hypothesis in (\ref{p2}) and
the fact that we have already dealt with the case in which two among
them are equal, we may assume for the moment that $1<c<b<a$. Let us
also denote $a\cdot b-c$ by $m$, $b\cdot c-a$ by $n$, and $c\cdot a
- b$ by $p$. Notice that $n<p<m$. By (\ref{a|b}), we must thus have
$n|p$, $n|m$, and $p|m$. Notice that $m-p = (b-c)\cdot (a+1)$ and
$m+p = (b+c)\cdot (a-1)$, so
\begin{equation} \label{b-ca}
p | (b-c)\cdot (a+1) \mbox{ and } p | (b+c)\cdot (a-1).
\end{equation}

One of $a+1$ and $a-1$ cannot be a multiple of $\overline{4}$, for
their difference is $\overline{2}$. If $a-1$ is not a multiple of
$\overline{4}$, then, since $p\cdot x= (b+c)\cdot (a-1)$ for some
$x>0$, and we have either $a-1=\overline{2}\cdot (\overline{2}\cdot
k +1)$ or $a-1=\overline{2}\cdot k +1$, we have $p\cdot x =
(b+c)\cdot \overline{2}\cdot (\overline{2}\cdot k +1)$ or $p\cdot x
= (b+c)\cdot (\overline{2}\cdot k +1)$. In both cases, by
(\ref{div2}), $\overline{2}\cdot k +1 | x$, i.e.,
$x=(\overline{2}\cdot k +1)\cdot y$, thus the two options are, after
canceling $\overline{2}\cdot k +1$ (by (\ref{canc})), $p\cdot y =
(b+c)\cdot \overline{2}$ or $p\cdot y = b+c$, thus in any case
$p\cdot y = \overline{2}\cdot(b+c)$ must hold for some $y$, and thus
\begin{equation} \label{p2b+c}
p\leq \overline{2}\cdot (b+c)
\end{equation}
If $a+1$ is not a multiple of $\overline{4}$, then we arrive
analogously to $p\cdot y = \overline{2}\cdot(b-c)$, and thus $p\leq
\overline{2}\cdot (b-c)$. So, in this case as well, (\ref{p2b+c})
holds.

Now, $b\cdot c + c= (b+1)\cdot c \leq a\cdot c = p+b\leq
\overline{2}\cdot(b+c)+b = \overline{3}\cdot b + \overline{2}\cdot
c$, thus $b\cdot c + c < \overline{3}\cdot b + \overline{2}\cdot c$,
thus, using A\ref{add}, $b\cdot c < \overline{3}\cdot b + c$, and,
given that $\overline{3}\cdot b + c <  \overline{4}\cdot b$, we get,
using A\ref{mul}, $c<\overline{4}$. Thus, we have only two
possibilities: (i) $c=\overline{2}$ and (ii)  $c=\overline{3}$.

Suppose (i) holds. Then we need to have $PT(a\cdot b -
\overline{2})$, $PT(\overline{2}\cdot a - b)$, and
$PT(\overline{2}\cdot b - a)$. If $a$ and $b$ were both even, then
$a\cdot b - \overline{2}$ would be a multiple of $\overline{2}$, but
not of $\overline{4}$, so we would need to have $a\cdot b -
\overline{2}=  \overline{2}$, which is impossible, since $b\geq
\overline{3}$ and $b\geq \overline{4}$. One can also easily notice
that $a$ and $b$ cannot both be odd, for else $a\cdot b -
\overline{2}$ would be odd, and thus would have to be $1$, which is
impossible for the reasons mentioned above. Thus the pair $(a,b)$
consists of an even and an odd number. Suppose $a$ were odd and $b$
were even, then $\overline{2}\cdot b - a$ would be odd, and thus
would have to be $1$. Thus $a=\overline{2}\cdot b-1$, and thus
$m=a\cdot b -c=\overline{2}\cdot b^2-b- \overline{2}$ and $p=c\cdot
a - b = \overline{3}\cdot b - \overline{2}$. Since $p|m$, we have
$\overline{3}\cdot b - \overline{2} | \overline{2}\cdot b^2-b-
\overline{2}$. Since
\begin{equation} \label{eq12}
\overline{9}\cdot(\overline{2}\cdot b^2-b- \overline{2}) =
(\overline{3}\cdot b - \overline{2})\cdot (\overline{6}\cdot b +1) -
\overline{16}
\end{equation}
we must have $\overline{3}\cdot b - \overline{2} | \overline{16}$.
Thus $\overline{3}\cdot b - \overline{2}\in \{1, \overline{2},
\overline{4},  \overline{8},  \overline{16}\}$. However, since
$b\geq \overline{3}$, we have $\overline{3}\cdot b -
\overline{2}\geq 7$, and thus we can have only $\overline{3}\cdot b
- \overline{2}= \overline{8}$, which has no solution $b$,  or
$\overline{3}\cdot b - \overline{2}= \overline{16}$, which means
$b=\overline{6}$ and $a=\overline{2}\cdot b-1= \overline{11}$. So,
in case $c=2$, we have only $(\overline{11}, \overline{6},
\overline{2})$ as solution.

Suppose now (ii) holds. Looking at (\ref{b-ca}) with $c=3$, we
notice that not both of $b-\overline{3}$ and $b+\overline{3}$ can be
multiples of $\overline{4}$ (given that their difference is
$\overline{6}$). If $4\nmid b-\overline{3}$, then $b-\overline{3}=
i\cdot (\overline{2}\cdot k +1)$ with $i\in \{1, \overline{2}\}$,
and (\ref{b-ca}) becomes $p\cdot x = i\cdot (\overline{2}\cdot k
+1)\cdot (a+1)$. By (\ref{div2}), $x=(\overline{2}\cdot k +1)\cdot
y$, so we have $p\cdot y = i\cdot (a+1)$, so $p\leq
\overline{2}\cdot (a+1)$. Similarly, if $4\nmid b+\overline{3}$,
then $p\cdot y = i\cdot (a-1)$, thus $p\leq  \overline{2}\cdot
(a-1)$. If  $4\nmid b-\overline{3}$, then we get $p\cdot y = i\cdot
(a+1)$, thus $p\leq  \overline{2}\cdot (a+1)$. So, in any case, we
have $p\leq  \overline{2}\cdot (a+1)$, i.e., $\overline{3}\cdot a -
b\leq  \overline{2}\cdot (a+1)$, which means $a-b\leq \overline{2}$.
Since we also have $1\leq a-b$, we can have only $a-b=1$ or $a-b=
\overline{2}$. If $a=b+1$, then $n=\overline{2}\cdot b -1$, which,
being odd and a power of $2$, must be $1$, which is not possible, as
it would imply $b=1$. If $a= b+\overline{2}$, then $m= (b-1)\cdot
(b+\overline{3})$, and thus we must have $PT(b-1)$ and
$PT(b+\overline{3})$. Since $(b+\overline{3})-(b-1)=\overline{4}$,
one of them must be $\overline{4}$, and, since $b\geq \overline{4}$,
that one cannot be $b+\overline{3}$, so it must be $b-1$, so
$b=\overline{5}$ and $a=\overline{7}$.

\end{proof}

\section{Pythagrean arithmetic is the right setting}

We may wonder whether we actually needed all of Pythagorean
Arithmetic to prove (\ref{p2}). From a methodological point of view,
we have argued that, in the absence of A\ref{rf}, the usual
properties of powers of $2$ would not hold, and thus the meaning of
the terms involved would be altered. In that sense Pythagorean
Arithmetic is the right theory in which the question regarding the
provability of (\ref{p2}) ought to be raised.

From a purely formal point of view, however, one is justified to ask
whether (\ref{p2}) does not follow from weaker assumptions. Our
proof already shows that it does. All we have used in it is PA$^-$,
A\ref{oe}, and (\ref{div2}). That this is less than what Pythagorean
Arithmetic asks can be seen by noticing that ${\mathcal C}({\mathbb
Q}(\sqrt{2})_{{\mathbb Z}} [X])$ is a model of PA$^-$, A\ref{oe},
and (\ref{div2}) (as there are no nonstandard powers of $2$ in it),
but not of Pythagorean Arithmetic (which is plain, as A\ref{rf}
fails for $m=X$ and $n=\sqrt{2}\cdot X$).

However, the weak theory of the odd and the even, PA$^- +$A\ref{oe},
is not strong enough to prove (\ref{p2}).
\begin{theorem}
PA$^-+$A\ref{oe} $\nvdash$ (\ref{p2}).
\end{theorem}
\begin{proof}

If $D$ is an ordered integral domain and $R$ is an ordered integral
domain containing $D$, then we denote by  $R_{D}[X, Y, Z]$ the ring
of polynomials in $X, Y,$ and $Z$, with free term in $D$ and with
all other coefficients in $R$, ordered by $\sum_{0\leq i,j,k \leq n}
c_{(i,j,k)}X^iY^jZ^k>0$ (here $c_{(0,0,0)}\in D$, and
$c_{(i,j,k)}\in K$ for all $1\leq i,j,k\leq n$) if and only if
$c_{(u,v,w)}>0$, where $(u,v,w)$ is the greatest element, in the
lexicographic ordering, among all the indexes $(i,j,k)$ of the
non-zero coefficients $c_{i,j,k}$ of the terms highest degree, i.\
e., for which $i+j+k$ is maximal (i.e., $(u,v,w)= \max \{(i,j,k)\,
:\, c_{(i,j,k)}\neq 0; i+j+k=d\}$, where $d$ is the degree of the
polynomial  $\sum_{0\leq i,j,k \leq n} c_{(i,j,k)}X^iY^jZ^k$ and
$\max$ is the greatest element in the lexicographic order). Let
${\mathcal C}(R_{D}[X, Y, Z])$ denote the positive cone of $R_{D}[X,
Y, Z]$.

Then ${\mathcal C}(R_{D}[X, Y, Z])$, with $R={\mathbb
Z}_{\frac{1}{2}}$ and $D={\mathbb Z}$, with $\left[\frac{\sum_{0\leq
i,j,k \leq n} c_{(i,j,k)}X^iY^jZ^k}{2}\right]=$\\
 $\sum_{0\leq i,j,k
\leq n, i+j+k\neq 0} \frac{c_{(i,j,k)}}{2}X^iY^jZ^k  +
\left[\frac{c_{(0,0,0)}}{2}\right]$, is a model of PA$^-+$A\ref{oe},
but not of (\ref{p2}), for all of $XY-Z$, $YZ-X$, $ZX-Y$ are
positive and are powers of two.

\end{proof}

Address: School of Mathematical and Natural Sciences (MC 2352)\\
Arizona State University - West campus\\
P. O. Box 37100\\
Phoenix, AZ 85306\\
U.S.A.\\
E-Mail: pamb@asu.edu
\end{document}